\documentclass[11pt,reqno]{amsart}

\usepackage[
left=0.5in,
right=0.5in,
top=0.7in,
bottom=0.7in
]{geometry}

\usepackage{amsmath,amsthm,amssymb,amscd,amsfonts,mathrsfs,color}
\usepackage[mathscr]{eucal}
\usepackage{mathtools}
\usepackage{pgfplots}
\usepackage{placeins}
\mathtoolsset{showonlyrefs}
\usepgfplotslibrary{fillbetween}
\pgfplotsset{compat=1.18}
\usepackage[bookmarksnumbered,colorlinks,plainpages]{hyperref}
\newtheorem{theorem}{Theorem}[section]
\newtheorem{proposition}[theorem]{Proposition}
\newtheorem{lemma}[theorem]{Lemma}
\newtheorem{corollary}[theorem]{Corollary}

\newtheorem{example}[theorem]{Example}
\theoremstyle{definition}

\numberwithin{equation}{section}

\newcommand{\BH}{\mathcal B(\mathcal H)}

\newcommand{\Real}{\operatorname{Re}}
\newcommand{\norm}[1]{\left\lVert #1\right\rVert}
\newcommand{\abs}[1]{\left\lvert #1\right\rvert}
\newcommand{\inner}[2]{\left\langle #1,#2\right\rangle}
\newcommand{\C}{\mathbb C}
\newcommand{\D}{\mathbb D}
\newcommand{\T}{\mathbb T}
\newcommand{\conv}{\operatorname{conv}}
\newcommand{\sgn}{\operatorname{sgn}}

\begin{document}

\title[Deformed $\rho$-teardrops and functional calculus]
{Deformed $\rho$-teardrops and a sharp
$f(0)$-dependent functional-calculus inequality}
\author[M. R. Lahwate and M. A. Ighachane]{Mohamed Rassan Lahwate$^{*}$ and Mohamed Amine Ighachane}

\address{Sciences and Technologies Team (ESTE),
			Higher School of Education and Training of El Jadida,
			Chouaib Doukkali University,
			El Jadida, Morocco.}
		\email{mohamedghassane94@gmail.com}

	\address{Sciences and Technologies Team (ESTE),
		Higher School of Education and Training of El Jadida,
		Chouaib Doukkali University,
		El Jadida, Morocco.}
	\email{mohamedamineighachane@gmail.com}

\subjclass[2020]{Primary 47A12, 47A20; Secondary 15A60, 47A30}
\keywords{$\rho$-operator radius, $\rho$-contraction, functional calculus,
Davis--Wielandt shell, numerical range, deformed teardrop}

\begin{abstract}
Let $1\leq\rho\leq2$. We determine the sharp pointwise bound for
$w_\rho(f(T))$ when $w_\rho(T)\leq1$, $f\in A(\overline\D)$,
$\norm{f}_{\infty,\D}\leq1$, and $\abs{f(0)}$ is prescribed.
Our approach is based on the exact universal Davis--Wielandt region of
$\mathcal C_\rho$ and its M\"obius deformation, which yields an exact
universal $\rho$-teardrop with a unique outer contact. This gives the
optimal functional-calculus inequality, recovers the sharp uniform
constant, and yields a rigidity characterization of equality.
\end{abstract}

\thanks{$^*$Corresponding author: Mohamed Rassan Lahwate}
	\maketitle
\section{Introduction and preliminaries}

The classes of operators admitting unitary $\rho$-dilations were introduced
by Sz.-Nagy and Foias~\cite{SzNagyFoias}. Let $\mathcal H$ be a complex
Hilbert space. For $\rho>0$, an operator $T\in\BH$ belongs to
$\mathcal C_\rho$ if there exist a Hilbert space
$\mathcal K\supseteq\mathcal H$ and a unitary
$U\in\mathcal B(\mathcal K)$ such that
\begin{equation}\label{eq:rho-dilation}
 T^n=\rho P_{\mathcal H}U^n\big|_{\mathcal H},
 \qquad n\geq1.
\end{equation}
The associated operator radius is defined by
\begin{equation}\label{eq:rho-radius}
 w_\rho(T)=\inf\{r>0:r^{-1}T\in\mathcal C_\rho\}.
\end{equation}
For $1\leq\rho\leq2$, $w_\rho$ is a norm, with
$w_1(T)=\norm{T}$ and $w_2(T)=w(T)$; see
Ando and Nishio~\cite{AndoNishio}. The structure and inequalities associated with $\rho$-operator radii
have continued to attract attention; see, among recent contributions,
\cite{KittanehZamani,BhuniaMoslehianZamani,GauWuConsequences}.

A natural functional-calculus problem is to determine the optimal bound
for $w_\rho(f(T))$ when $w_\rho(T)\leq1$ and $f$ belongs to the unit
ball of the disk algebra $A(\overline{\D})$. At the numerical-radius
endpoint $\rho=2$, the normalized case $f(0)=0$ is covered by the
classical Berger--Stampfli mapping theorem~\cite{BergerStampfli}; see
also~\cite{KlajaMashreghiRansford}. For $\rho$-contractions, the
corresponding normalized mapping property is known for rational
functions with poles off $\overline{\D}$; see
\cite{SzNagyFoias,Williams,SchwenningerDeVries}. Using the spectral-set
closure principle of Badea et al.~\cite{BadeaCrouzeixKlaja}, we then
pass to the disk algebra $A(\overline{\D})$ by uniform approximation.

When no condition is imposed on $f(0)$, Badea et
al.~\cite{BadeaCrouzeixKlaja} obtained the sharp uniform estimate
\begin{equation}\label{eq:uniform-intro}
 w_\rho(f(T))
 \leq
 C_\rho
 :=
 \frac{\rho^2+1+(\rho-1)\sqrt{5\rho^2+2\rho+1}}
      {2\rho^2}.
\end{equation}

At the endpoint $\rho=2$, substantially more information is available:
the teardrop theorem gives the pointwise estimate
\begin{equation}\label{eq:Drury-pointwise}
 w(f(T))
 \leq
 1+s-s^2,
 \qquad
 s=\abs{f(0)}.
\end{equation}
See~\cite{Drury}. Our purpose is to extend the pointwise estimate
\eqref{eq:Drury-pointwise} to the full range
$1\leq\rho\leq2$ by determining the optimal bound in terms of the
single parameter $\abs{f(0)}$. More precisely, we prove
\begin{equation}\label{eq:main-intro}
 w_\rho(f(T))
 \leq
 \ell_\rho\bigl(\abs{f(0)}\bigr),
 \qquad
 1\leq\rho\leq2.
\end{equation}

For $0\leq s\leq1$, the function $\ell_\rho$ is given by
\begin{equation}\label{eq:B-ell}
 B_\rho(s)
 =
 \rho(1-s^2)+2(\rho-1)s,
 \qquad
 \ell_\rho(s)
 =
 \frac{
 B_\rho(s)+
 \sqrt{B_\rho(s)^2+4\rho(2-\rho)s^2}
 }
 {2\rho}.
\end{equation}
Equivalently, $t=\ell_\rho(s)$ is the positive root of
\begin{equation}\label{eq:ell-quadratic}
 \rho t^2-B_\rho(s)t-(2-\rho)s^2=0.
\end{equation}
The proof of~\eqref{eq:main-intro} reduces the disk-algebra problem to
the M\"obius family
\[
 g_s(z)=\frac{s+z}{1+sz}.
\]
For $C\in\mathcal C_\rho$, we determine the exact universal
Davis--Wielandt region, compute its M\"obius deformation, and apply the
deformed numerical-range map. The resulting universal $\rho$-teardrop has outer radius $\ell_\rho(s)$
and a unique outer contact. These two facts yield the sharp bound and
the equality characterization.

The estimate~\eqref{eq:main-intro} is sharp for every prescribed value
of $\abs{f(0)}$.

Section~2 determines the universal Davis--Wielandt region and its support
geometry. Section~3 constructs the deformed $\rho$-teardrop and determines
its outer radius and unique contact, together with an ordinary
numerical-range companion. Section~4 derives~\eqref{eq:main-intro} and
recovers the sharp uniform constant~\eqref{eq:uniform-intro}.
Section 5 characterizes equality through the unique extremal
Davis--Wielandt contact and derives the corresponding criterion for
the original pair $(T,f)$.

Throughout the paper, we write $p=\rho-1$ and $q=2-\rho$.
Whenever division by $p$ or $q$ occurs, the corresponding endpoint
$\rho=1$ or $\rho=2$ is treated separately. We shall use the following two
positivity criteria. The boundary criterion goes back to Sz.-Nagy and
Foias~\cite{SzNagyFoias}, while the disk form is
\cite[Theorem~2.2]{SchwenningerDeVries}.

\begin{lemma}\label{lem:positivity}
Let $A\in\BH$ and $r>0$.

\noindent{\rm (i)} One has $w_\rho(A)\leq r$ if and only if
\begin{equation}\label{eq:boundary-pencil}
 Q_A(\zeta,r)
 :=
 \rho r^2I-2pr\Real(\overline\zeta A)-qA^*A
 \geq0,
 \qquad \zeta\in\T.
\end{equation}

\noindent{\rm (ii)} An operator $C$ belongs to $\mathcal C_\rho$ if and
only if
\begin{equation}\label{eq:disk-pencil}
 P_C(u)
 :=
 I-\frac{2p}{\rho}\Real(\overline u C)
 -\frac{1-p^2\abs{u}^2}{\rho^2}C^*C
 \geq0,
 \qquad u\in\overline\D.
\end{equation}
\end{lemma}

\begin{proof}
For part~{\rm (i)}, the boundary characterization applied to $r^{-1}A$
gives
\[
 \rho I-\frac{2p}{r}\Real(\overline\zeta A)
 -\frac{q}{r^2}A^*A\geq0,
 \qquad \zeta\in\T.
\]
Multiplication by $r^2$ gives~\eqref{eq:boundary-pencil}, and the
converse follows by reversing the same calculation.

For part~{\rm (ii)}, the cited disk test applies first when
$\sigma(C)\subset\D$. Let $C\in\mathcal C_\rho$ and $0<v<1$. Since
\[
 v\bigl(\rho I-2p\Real(\overline\zeta C)-qC^*C\bigr)
 +\rho(1-v)I+qv(1-v)C^*C\geq0,
\]
the boundary criterion shows that $vC\in\mathcal C_\rho$.
Moreover,~\eqref{eq:rho-dilation} gives $\norm{C^n}\leq\rho$, and hence
$\sigma(vC)\subset\D$. Applying the disk test to $vC$ and letting
$v\uparrow1$ yields~\eqref{eq:disk-pencil}.

Conversely, by continuity,~\eqref{eq:disk-pencil} extends to
$\abs{u}=1$. Multiplying by $\rho$ then gives
~\eqref{eq:boundary-pencil} with $r=1$.
\end{proof}

We shall use the factorization from~\cite[Theorem~2.2]{SchwenningerDeVries}.
Assume $1<\rho\leq2$. For $C\in\mathcal C_\rho$ and $w\in\D$, one has
\[
 \left|\frac{1+pw}{\rho}\right|
 \leq
 \frac{1+p\abs{w}}{\rho}
 <1.
\]
Since the spectral radius of $C$ is at most one, it follows that
\[
 S_C(w)
 =
 \left(I-\frac{1+pw}{\rho}C\right)^{-1}
\]
is well defined. Set
\[
 F_C(w)
 =
 I+\frac{2}{\rho}CS_C(w).
\]
Then
\begin{equation}\label{eq:positive-real-factorization}
 \Real F_C(w)
 =
 S_C(w)^*
 \left[
 I-\frac{2p}{\rho}\Real(wC)
 -\frac{1-p^2\abs{w}^2}{\rho^2}C^*C
 \right]
 S_C(w)
 \geq0.
\end{equation}

Let $A(\overline\D)$ denote the disk algebra, namely the algebra of
functions continuous on $\overline\D$ and holomorphic in $\D$.
Polynomials are uniformly dense in $A(\overline\D)$, and hence
\[
 R(\overline\D)=A(\overline\D),
\]
where $R(\overline\D)$ is the uniform closure of the rational functions
with poles off $\overline\D$.

For $C\in\mathcal C_\rho$, the rational functional calculus extends
continuously to $A(\overline\D)$. Indeed,
\cite[Theorem~5.2]{BadeaCrouzeixKlaja} gives
\[
 w_\rho(r(C))
 \leq
 C_\rho\norm{r}_{\infty,\D}
\]
for every rational function $r$ with poles off $\overline\D$. Since
$\norm{A}\leq\rho w_\rho(A)$, it follows that
\begin{equation}\label{eq:bounded-disk-calculus}
 \norm{r(C)}
 \leq
 \rho C_\rho\norm{r}_{\infty,\D}
\end{equation}
for every such $r$.
By~\eqref{eq:bounded-disk-calculus}, the rational calculus has a unique continuous
extension to $A(\overline{\D})$. This extension is a unital algebra
homomorphism, and we use $h(C)$ to denote its value at
$h\in A(\overline{\D})$.

We shall use the following normalized mapping theorem for
$\rho$-contractions. The rational case follows from
\cite{SzNagyFoias,Williams}; see also
\cite[Theorem~3.1]{SchwenningerDeVries}. The disk-algebra formulation
below follows by uniform approximation, using the spectral-set closure
principle of \cite{BadeaCrouzeixKlaja}.

\begin{theorem}\label{thm:zero-mapping}
Let $1\leq\rho\leq2$, let $C\in\mathcal C_\rho$, and let
$h\in A(\overline{\D})$. If
\[
 \norm{h}_{\infty,\D}\leq1,
 \qquad
 h(0)=0,
\]
then
\[
 h(C)\in\mathcal C_\rho.
\]
\end{theorem}

\begin{proof}
Choose polynomials $p_n$ converging uniformly to $h$ on
$\overline\D$, and put
\[
 \widetilde p_n=p_n-p_n(0),
 \qquad
 a_n=\max\left\{1,\norm{\widetilde p_n}_{\infty,\D}\right\},
 \qquad
 q_n=\frac{\widetilde p_n}{a_n}.
\]
Since $p_n(0)\to h(0)=0$ and $\widetilde p_n\to h$ uniformly on
$\overline\D$, one has $a_n\to1$ and
\[
 q_n\longrightarrow h
 \qquad\text{uniformly on }\overline\D.
\]
Moreover,
\[
 q_n(0)=0,
 \qquad
 \norm{q_n}_{\infty,\D}\leq1.
\]
The rational normalized mapping theorem gives
$q_n(C)\in\mathcal C_\rho$ for every $n$. By the continuous
disk-algebra functional calculus defined above,
\[
 q_n(C)\longrightarrow h(C)
 \qquad\text{in operator norm}.
\]
Finally, $\mathcal C_\rho$ is norm closed: this follows directly from
the boundary inequalities~\eqref{eq:boundary-pencil}. Hence
$h(C)\in\mathcal C_\rho$.
\end{proof}

We next recall the part of the deformed numerical-range theory needed
below. This range was introduced by Pagacz et al.\
\cite{PagaczPietrzyckiWojtylak}. For a unit vector $x$, write
$z=\inner{Ax}{x}$ and $b=\norm{Ax}^2$. For $1\leq\rho\leq2$, define
\begin{equation}\label{eq:Phi-rho}
 \Phi_\rho(z,b)
 =
 \begin{cases}
 \displaystyle
 \frac1\rho
 \left(
 pz+\frac{z}{\abs{z}}
 \sqrt{p^2\abs{z}^2+\rho qb}
 \right),
 & z\neq0,\\[6pt]
 0,
 & z=0.
 \end{cases}
\end{equation}
The second line is not the continuous extension of the first; it is the
value prescribed by the original definition when $\inner{Ax}{x}=0$.
The set $W^\rho(A)$ is the closed convex hull of these generating points.
By~\cite[Theorem~9 and Corollary~10]{PagaczPietrzyckiWojtylak},
\begin{equation}\label{eq:deformed-facts}
 W^2(A)=\overline{W(A)},
 \qquad
 A\in\mathcal C_\rho
 \Longleftrightarrow
 W^\rho(A)\subseteq\overline\D.
\end{equation}
The same work also gives the homogeneity
$W^\rho(\alpha A)=\alpha W^\rho(A)$. Consequently,
\eqref{eq:rho-radius} and~\eqref{eq:deformed-facts} yield
\begin{equation}\label{eq:gauge-Wrho}
 w_\rho(A)
 =
 \max\{\abs{\lambda}:\lambda\in W^\rho(A)\}.
\end{equation}
Indeed, after scaling by $r>0$, the right-hand side is at most one
precisely when $r^{-1}A\in\mathcal C_\rho$; taking the infimum over $r$
gives the identity. For $z\neq0$, the modulus of~\eqref{eq:Phi-rho} is
the positive root $R$ of
\begin{equation}\label{eq:Phi-quadratic}
 \rho R^2-2p\abs{z}R-qb=0.
\end{equation}

\section{The universal Davis--Wielandt region}

Finite-dimensional and operator-matrix aspects of $\rho$-contractions
have been studied in \cite{OkuboSpitkovsky,CassierZerouali}. Our aim in
this section is to determine the exact Davis--Wielandt region generated
by the whole class $\mathcal C_\rho$, together with the support geometry
that will be needed later.

For $1<\rho<2$, define
\[
 a_\rho=\frac{\rho p}{1+p^2},
 \qquad
 \beta_\rho(a)
 =
 \begin{cases}
 \displaystyle
 \frac{\rho}{2}
 \left(\rho+\sqrt{\rho^2-4a^2}\right),
 & 0\leq a\leq a_\rho,\\[6pt]
 \displaystyle
 \frac{\rho-2pa}{q},
 & a_\rho\leq a\leq1.
 \end{cases}
\]
For $\rho=2$, set
\[
 \beta_2(a)=2+2\sqrt{1-a^2},
 \qquad 0\leq a\leq1.
\]

For an operator $A$, define its Davis--Wielandt shell by
\[
 DW(A)
 =
 \{(\inner{Ax}{x},\norm{Ax}^2):\norm{x}=1\}.
\]
Following Davis~\cite{Davis}, we regard $DW(A)$ in complex--real
coordinates. Throughout the universal constructions
below, the underlying complex Hilbert space is allowed to vary.
Accordingly, let $\mathfrak D_\rho$ denote the union of
$\operatorname{DW}(C)$ over all $C\in\mathcal C_\rho$ and all such
Hilbert spaces.

We now determine this universal region exactly.

\begin{theorem}\label{thm:exact-DW}
For $1<\rho\leq2$,
\[
 \mathfrak D_\rho
 =
 \left\{
 (z,b)\in\C\times\mathbb R:
 \abs{z}\leq1,\;
 \abs{z}^2\leq b\leq\beta_\rho(\abs{z})
 \right\}.
\]
Moreover, every point of $\mathfrak D_\rho$ is realized by an operator
acting on a Hilbert space of dimension at most two.
\end{theorem}

Before proving the theorem, we record the compression property needed in
the realization argument.

\begin{lemma}\label{lem:compression}
If $C\in\mathcal C_\rho$ and $\mathcal M\subseteq\mathcal H$ is a closed
subspace, then
\[
 P_{\mathcal M}C\big|_{\mathcal M}\in\mathcal C_\rho.
\]
\end{lemma}

\begin{proof}
Set
\[
 B=P_{\mathcal M}C\big|_{\mathcal M}.
\]
For $x\in\mathcal M$ and $\zeta\in\T$, we have
$\inner{Bx}{x}=\inner{Cx}{x}$ and
$\norm{Bx}=\norm{P_{\mathcal M}Cx}\leq\norm{Cx}$. Since $q\geq0$,
\begin{align*}
 \inner{Q_B(\zeta,1)x}{x}
 &=
 \rho\norm{x}^2
 -2p\Real\bigl(\overline\zeta\inner{Cx}{x}\bigr)
 -q\norm{P_{\mathcal M}Cx}^2\\
 &\geq
 \rho\norm{x}^2
 -2p\Real\bigl(\overline\zeta\inner{Cx}{x}\bigr)
 -q\norm{Cx}^2\\
 &=
 \inner{Q_C(\zeta,1)x}{x}
 \geq0.
\end{align*}
Hence $Q_B(\zeta,1)\geq0$ for every $\zeta\in\T$, and
Lemma~\ref{lem:positivity}{\rm(i)} gives $B\in\mathcal C_\rho$.
\end{proof}

\begin{proof}[Proof of Theorem~\ref{thm:exact-DW}]
We first prove the inclusion
\[
 \mathfrak D_\rho
 \subseteq
 \left\{
 (z,b)\in\C\times\mathbb R:
 \abs{z}\leq1,\;
 \abs{z}^2\leq b\leq\beta_\rho(\abs{z})
 \right\}.
\]
Let $C\in\mathcal C_\rho$, let $x$ be a unit vector, and put
\[
 a=\abs{\inner{Cx}{x}},
 \qquad
 b=\norm{Cx}^2.
\]
By Cauchy--Schwarz, $a^2\leq b$. Choose $u\in\overline\D$ with
$\abs{u}=v$ and phase such that
$\Real(\overline u\inner{Cx}{x})=va$. Applying
\eqref{eq:disk-pencil} to $x$ gives
\begin{equation}\label{eq:scalar-disk}
 1-\frac{b}{\rho^2}
 -\frac{2p}{\rho}va
 +\frac{p^2}{\rho^2}v^2b
 \geq0,
 \qquad 0\leq v\leq1.
\end{equation}
Taking $v=1$ and using $b\geq a^2$ yields
$qa^2+2pa\leq\rho=q+2p$, hence $a\leq1$.

Assume first that $1<\rho<2$ and $b>0$. The left-hand side of
\eqref{eq:scalar-disk} attains its minimum on $[0,1]$ at
\[
 v_0=\min\left\{1,\frac{\rho a}{pb}\right\}.
\]
If $\rho a\leq pb$, substitution of the interior critical point gives
\[
 \frac{b}{\rho^2}+\frac{a^2}{b}\leq1,
\]
and therefore
\[
 b
 \leq
 \frac{\rho}{2}
 \left(\rho+\sqrt{\rho^2-4a^2}\right).
\]
The case condition is compatible with this estimate only when
$a\leq a_\rho$. Indeed,
\[
 \frac{\rho a}{p}
 \leq
 \frac{\rho}{2}
 \left(\rho+\sqrt{\rho^2-4a^2}\right)
 \quad\Longleftrightarrow\quad
 a\leq\frac{\rho p}{1+p^2}=a_\rho.
\]

If $\rho a\geq pb$, evaluating~\eqref{eq:scalar-disk} at $v=1$ gives
\[
 b\leq\frac{\rho-2pa}{q}.
\]
For $a\geq a_\rho$, this is the affine branch of $\beta_\rho$. When
$a\leq a_\rho$, the case condition gives
$b\leq\rho a/p$, and the preceding equivalence shows that this quantity
is bounded above by the curved branch of $\beta_\rho$. Hence
$b\leq\beta_\rho(a)$ in every case.

For $\rho=2$, the same minimization argument gives
$b/4+a^2/b\leq1$ in the interior case, hence
\[
 b\leq2+2\sqrt{1-a^2}.
\]
If the minimum occurs at $v=1$, then
$b\leq2a\leq2+2\sqrt{1-a^2}$. This proves the first inclusion.

We now prove the reverse inclusion. We refer to
$
 \bigl(\inner{Ax}{x},\norm{Ax}^2\bigr)
$
as the Davis--Wielandt pair of $A$ at $x$. Consider the canonical
nilpotent
\[
 N_\rho=
 \begin{pmatrix}
 0&\rho\\
 0&0
 \end{pmatrix}.
\]
For $\zeta\in\T$,
\begin{equation}\label{eq:N-pencil}
 Q_{N_\rho}(\zeta,1)
 =
 \rho
 \begin{pmatrix}
 1&-p\overline\zeta\\
 -p\zeta&p^2
 \end{pmatrix}
 \geq0,
\end{equation}
so $N_\rho\in\mathcal C_\rho$ by
Lemma~\ref{lem:positivity}{\rm(i)}.

We first realize the curved part of the upper boundary. For
$0\leq a\leq a_\rho$ when $\rho<2$, and for $0\leq a\leq1$ when
$\rho=2$, set
\[
 u=\frac12\left(1+\sqrt{1-\frac{4a^2}{\rho^2}}\right),
 \qquad
 x=(\sqrt{1-u},\sqrt u)^{\mathsf T}.
\]
The Davis--Wielandt pair of $N_\rho$ at $x$ is
\[
 \left(
 a,\,
 \frac{\rho}{2}
 \left(\rho+\sqrt{\rho^2-4a^2}\right)
 \right),
\]
and hence $N_\rho$ realizes the curved upper boundary.

Suppose now that $1<\rho<2$. At the transition point, put
\[
 b_\rho=\frac{\rho^2}{1+p^2},
 \qquad
 x_\rho=\frac1{\sqrt{1+p^2}}(p,1)^{\mathsf T}.
\]
The Davis--Wielandt pair of $N_\rho$ at $x_\rho$ is
$(a_\rho,b_\rho)$. For $a_\rho\leq a\leq1$, let
\[
 \lambda=\frac{1-a}{1-a_\rho},
 \qquad
 y=\sqrt{\lambda}\,x_\rho\oplus\sqrt{1-\lambda}.
\]
Since finite direct sums of operators in $\mathcal C_\rho$ remain in
$\mathcal C_\rho$, we may consider $N_\rho\oplus I$. Its
Davis--Wielandt pair at $y$ is
\[
 \left(a,\lambda b_\rho+1-\lambda\right)
 =
 \left(a,\beta_\rho(a)\right).
\]
Thus every point of the upper boundary is realized by an operator in
$\mathcal C_\rho$.

It remains to fill the region below this boundary. Fix
$a^2\leq b\leq\beta_\rho(a)$ and let $(C_a,x_a)$ realize the upper
boundary pair $(a,\beta_\rho(a))$. The lower endpoint $(a,a^2)$ is
realized by the scalar operator $aI$. If
$a^2<b<\beta_\rho(a)$, set
\[
 \mu=\frac{b-a^2}{\beta_\rho(a)-a^2},
 \qquad
 v=\sqrt{\mu}\,x_a\oplus\sqrt{1-\mu}.
\]
Then the Davis--Wielandt pair of $C_a\oplus aI$ at $v$ is $(a,b)$.
Hence every pair with
$a^2\leq b\leq\beta_\rho(a)$ is realized by an operator in
$\mathcal C_\rho$.

For $z\neq0$, multiplication of the realizing operator by the
unimodular scalar $z/\abs{z}$ changes its first Davis--Wielandt
coordinate from $\abs{z}$ to $z$ and leaves the second unchanged. The
case $z=0$ is already covered by taking $a=0$. This proves the reverse
inclusion.

It remains only to verify the dimension assertion. Let $C\in\mathcal
C_\rho$ and a unit vector $x$ realize any of the pairs constructed
above, and set
$\mathcal M=\operatorname{span}\{x,Cx\}$. Then $\dim\mathcal M\leq2$, and for
$B=P_{\mathcal M}C|_{\mathcal M}$ one has $Bx=Cx$. Thus $B$ has the
same Davis--Wielandt pair at $x$ as $C$. By
Lemma~\ref{lem:compression}, $B\in\mathcal C_\rho$. Consequently every
point of $\mathfrak D_\rho$ is realized in dimension at most two.
\end{proof}

The endpoint $\rho=1$ has the analogous region
$|z|\leq1$, $|z|^2\leq b\leq1$. Its upper boundary follows from the
unitary $\operatorname{diag}(1,-1)$, and the preceding direct-sum and
compression argument gives all intermediate points.

The following support formula will locate the contact point used in
Sections~3 and~5.

\begin{lemma}\label{lem:beta-support}
Let $1<\rho\leq2$. For $u\geq0$ and $k>0$,
\begin{equation}\label{eq:beta-support}
 \max_{0\leq a\leq1}\{2ua+k\beta_\rho(a)\}
 =
 \begin{cases}
 \displaystyle
 \frac{\rho}{2}
 \left(
 \rho k+\sqrt{\rho^2k^2+4u^2}
 \right),
 & qu\leq pk,\\[6pt]
 2u+k,
 & qu\geq pk.
 \end{cases}
\end{equation}
For $\rho=2$, only the first alternative occurs. If $qu<pk$, the
maximizer is unique and equals
\begin{equation}\label{eq:beta-maximizer}
 a_*=
 \frac{\rho u}{\sqrt{4u^2+\rho^2k^2}}.
\end{equation}
Moreover, $a_*<a_\rho$ when $1<\rho<2$, whereas $a_*<1$ when
$\rho=2$.
\end{lemma}

\begin{proof}
Put $F(a)=2ua+k\beta_\rho(a)$. Assume first that $1<\rho<2$. On the
curved branch $0\leq a\leq a_\rho$, write
\[
 F_{\rm c}(a)
 =
 2ua+\frac{\rho k}{2}
 \left(\rho+\sqrt{\rho^2-4a^2}\right).
\]
For $0\leq a<a_\rho$, differentiation gives
\[
 F_{\rm c}'(a)
 =
 2u-\frac{2\rho ka}{\sqrt{\rho^2-4a^2}},
 \qquad
 F_{\rm c}''(a)
 =
 -\frac{2\rho^3k}{(\rho^2-4a^2)^{3/2}}<0.
\]
Hence $F_{\rm c}$ is strictly concave. Its critical point is
\[
 a_*=\frac{\rho u}{\sqrt{4u^2+\rho^2k^2}}.
\]
At this point,
\[
 \sqrt{\rho^2-4a_*^2}
 =
 \frac{\rho^2k}{\sqrt{4u^2+\rho^2k^2}},
\]
and therefore
\begin{equation}\label{eq:curved-support-value}
 F_{\rm c}(a_*)
 =
 \frac\rho2
 \left(
 \rho k+\sqrt{\rho^2k^2+4u^2}
 \right).
\end{equation}

Since $a_\rho=\rho p/(1+p^2)$ and $1-p^2=\rho q$, direct squaring gives
\begin{equation}\label{eq:critical-transition}
 a_*\leq a_\rho
 \quad\Longleftrightarrow\quad
 q^2u^2\leq p^2k^2
 \quad\Longleftrightarrow\quad
 qu\leq pk.
\end{equation}

On the affine branch $a_\rho\leq a\leq1$, one has
\[
 F_{\rm a}(a)
 =
 2ua+k\frac{\rho-2pa}{q},
 \qquad
 F_{\rm a}'(a)
 =
 2\left(u-\frac{pk}{q}\right).
\]
The two branch formulas agree at $a=a_\rho$.

If $qu<pk$, then $a_*<a_\rho$ by~\eqref{eq:critical-transition}, while
\[
 F_{\rm a}'(a)
 =
 2\left(u-\frac{pk}{q}\right)<0.
\]
Hence $F_{\rm a}$ is strictly decreasing on $[a_\rho,1]$. Since
$F_{\rm c}$ is strictly concave and attains its unique maximum at
$a_*$, the unique maximizer of $F$ is $a_*$, and
\eqref{eq:curved-support-value} gives the first value in
\eqref{eq:beta-support}.

If $qu>pk$, then $a_*>a_\rho$. Since $F_{\rm c}'$ vanishes only at
$a_*$, it follows that $F_{\rm c}'(a)>0$ for
$0\leq a\leq a_\rho$. Moreover,
\[
 F_{\rm a}'(a)
 =
 2\left(u-\frac{pk}{q}\right)>0.
\]
Thus $F$ is increasing on both branches and attains its maximum at
$a=1$. Since $\beta_\rho(1)=1$, this maximum is $2u+k$.

Finally, if $qu=pk$, then $a_*=a_\rho$ and
$F_{\rm a}'(a)=0$. Hence $F_{\rm a}$ is constant on
$[a_\rho,1]$, and the two expressions in
\eqref{eq:beta-support} coincide. This proves
\eqref{eq:beta-support} for $1<\rho<2$, together with the asserted
uniqueness when $qu<pk$.

Let now $\rho=2$. Then
\[
 F(a)=2ua+2k+2k\sqrt{1-a^2},
 \qquad 0\leq a\leq1.
\]
For $0\leq a<1$,
\[
 F'(a)
 =
 2u-\frac{2ka}{\sqrt{1-a^2}},
 \qquad
 F''(a)
 =
 -\frac{2k}{(1-a^2)^{3/2}}<0.
\]
Thus $F$ is strictly concave, and its unique critical point is
\[
 a_*=\frac{u}{\sqrt{u^2+k^2}}<1.
\]
Substitution gives
\[
 F(a_*)
 =
 2k+2\sqrt{u^2+k^2}.
\]
This is exactly the first expression in~\eqref{eq:beta-support} when
$\rho=2$, and it agrees with~\eqref{eq:beta-maximizer}.
\end{proof}

\section{M\"obius deformation and the deformed $\rho$-teardrop}

For $0\leq s<1$, let $g_s$ denote the M\"obius automorphism of $\D$
defined by
\[
 g_s(z)=\frac{s+z}{1+sz}.
\]

\subsection{Exact Davis--Wielandt deformation}

Let $C\in\mathcal C_\rho$, and put $X=g_s(C)$ and $D=I+sC$.
Since the spectral radius of $C$ is at most one, $D$ is invertible. Define
\[
 \Psi_s(z,b)
 =
 \left(
 \frac{s(1+b)+z+s^2\overline z}{1+2s\Real z+s^2b},
 \frac{s^2+2s\Real z+b}{1+2s\Real z+s^2b}
 \right).
\]

\begin{proposition}\label{prop:exact-deformation}
For every $C\in\mathcal C_\rho$,
\begin{equation}\label{eq:exact-deformation}
 DW(g_s(C))=\Psi_s(DW(C)).
\end{equation}
Consequently,
\begin{equation}\label{eq:universal-deformation}
 \bigcup_{C\in\mathcal C_\rho}DW(g_s(C))
 =
 \Psi_s(\mathfrak D_\rho),
\end{equation}
and every point on the right is realized on a Hilbert space of dimension
at most two.
\end{proposition}

\begin{proof}
Let $y$ be a unit vector, set $z=\inner{Cy}{y}$ and
$b=\norm{Cy}^2$, and define $x=Dy/\norm{Dy}$. Since
$d:=\norm{Dy}^2=1+2s\Real z+s^2b$, one obtains
\begin{align*}
 \inner{Xx}{x}
 &=
 \frac{\inner{(sI+C)y}{Dy}}{d}
 =
 \frac{s(1+b)+z+s^2\overline z}{d},\\
 \norm{Xx}^2
 &=
 \frac{s^2+2s\Real z+b}{d}.
\end{align*}
Thus the Davis--Wielandt pair of $X$ at $x$ is
$\Psi_s(z,b)$, which proves
\[
 DW(g_s(C))\subseteq\Psi_s(DW(C)).
\]

Conversely, let $x$ be a unit vector and normalize $D^{-1}x$ to obtain
a unit vector $y$. Then $x=Dy/\norm{Dy}$ up to a unimodular scalar, and
the same calculation shows that the Davis--Wielandt pair of $X$ at $x$
belongs to $\Psi_s(DW(C))$. Hence~\eqref{eq:exact-deformation} holds.

Taking the union over $C\in\mathcal C_\rho$ and using
Theorem~\ref{thm:exact-DW} gives~\eqref{eq:universal-deformation}.
The same theorem also gives the two-dimensional realization.
\end{proof}

\subsection{The universal deformed teardrop}

If $(Z,B)=\Psi_s(z,b)$, define
$\Theta_{\rho,s}(z,b)=\Phi_\rho(Z,B)$, using the second line of
\eqref{eq:Phi-rho} when $Z=0$. Set
\[
 \mathcal T_{\rho,s}
 =
 \overline{\conv}\bigl(\Theta_{\rho,s}(\mathfrak D_\rho)\bigr).
\]

The next result identifies this set as the exact universal deformed
numerical range associated with the Möbius map $g_s$.

\begin{theorem}\label{thm:universal-teardrop}
For $1<\rho\leq2$ and $0\leq s<1$,
\[
 \overline{\bigcup_{C\in\mathcal C_\rho}W^\rho(g_s(C))}
 =
 \mathcal T_{\rho,s}.
\]
\end{theorem}

\begin{proof}
For $A\in\BH$, denote by
\[
 \mathcal G_\rho(A)
 =
 \left\{
 \Phi_\rho\bigl(\inner{Ax}{x},\norm{Ax}^2\bigr):
 \norm{x}=1
 \right\}
\]
the set of generating points of $W^\rho(A)$. By definition,
\begin{equation}\label{eq:deformed-generating-hull}
 W^\rho(A)
 =
 \overline{\conv}\mathcal G_\rho(A).
\end{equation}

We first prove the inclusion from left to right. Let
$C\in\mathcal C_\rho$. Proposition~\ref{prop:exact-deformation} gives
\begin{equation}\label{eq:generator-deformation}
 \mathcal G_\rho(g_s(C))
 =
 \Theta_{\rho,s}(DW(C)).
\end{equation}
Since $DW(C)\subseteq\mathfrak D_\rho$, it follows from
\eqref{eq:deformed-generating-hull} that
\[
 W^\rho(g_s(C))
 \subseteq
 \overline{\conv}\Theta_{\rho,s}(\mathfrak D_\rho)
 =
 \mathcal T_{\rho,s}.
\]
Taking the union over $C\in\mathcal C_\rho$ and then the closure yields
\[
 \overline{\bigcup_{C\in\mathcal C_\rho}W^\rho(g_s(C))}
 \subseteq
 \mathcal T_{\rho,s}.
\]

For the reverse inclusion, consider a finite convex combination
\begin{equation}\label{eq:finite-theta-combination}
 \theta
 =
 \sum_{j=1}^n\alpha_j\theta_j,
 \qquad
 \theta_j=\Theta_{\rho,s}(z_j,b_j),
 \qquad
 \alpha_j\geq0,
 \qquad
 \sum_{j=1}^n\alpha_j=1,
\end{equation}
where $(z_j,b_j)\in\mathfrak D_\rho$. By
Theorem~\ref{thm:exact-DW}, for each $j$ there exist a Hilbert space
$\mathcal H_j$, an operator $C_j\in\mathcal C_\rho$, and a unit vector
$y_j\in\mathcal H_j$ whose Davis--Wielandt pair is $(z_j,b_j)$.

Put $X_j=g_s(C_j)$ and
\[
 x_j
 =
 \frac{(I+sC_j)y_j}{\norm{(I+sC_j)y_j}}.
\]
Proposition~\ref{prop:exact-deformation} shows that the
Davis--Wielandt pair of $X_j$ at $x_j$ is
$\Psi_s(z_j,b_j)$. Therefore
\begin{equation}\label{eq:theta-is-generator}
 \theta_j\in\mathcal G_\rho(X_j).
\end{equation}

Now set
\[
 \mathcal H=\bigoplus_{j=1}^n\mathcal H_j,
 \qquad
 C=\bigoplus_{j=1}^n C_j.
\]
For every $\zeta\in\T$,
\[
 Q_C(\zeta,1)
 =
 \bigoplus_{j=1}^nQ_{C_j}(\zeta,1)
 \geq0.
\]
Hence Lemma~\ref{lem:positivity}{\rm(i)} gives
$C\in\mathcal C_\rho$. Moreover, finite direct sums commute with the
rational functional calculus, so
\[
 g_s(C)
 =
 \bigoplus_{j=1}^nX_j.
\]

Regard $x_j$ as a unit vector in the $j$th summand of $\mathcal H$.
Its Davis--Wielandt pair for $g_s(C)$ is the same as for $X_j$, and
hence~\eqref{eq:theta-is-generator} gives
$\theta_j\in\mathcal G_\rho(g_s(C))$ for every $j$. Consequently,
\eqref{eq:finite-theta-combination} and
\eqref{eq:deformed-generating-hull} imply
\[
 \theta
 \in
 \conv\mathcal G_\rho(g_s(C))
 \subseteq
 W^\rho(g_s(C)).
\]
Thus
\[
 \conv\Theta_{\rho,s}(\mathfrak D_\rho)
 \subseteq
 \bigcup_{C\in\mathcal C_\rho}W^\rho(g_s(C)).
\]
Taking closures gives the reverse inclusion.
\end{proof}

Theorem~\ref{thm:universal-teardrop} identifies the exact universal
deformed numerical range. We now determine its outer radius and show
that the corresponding outer circle has a unique point of contact.

\begin{theorem}\label{thm:radial-teardrop}
Let $1<\rho\leq2$ and $0<s<1$. Then
\begin{equation}\label{eq:teardrop-radius}
 \max\{\abs{\lambda}:\lambda\in\mathcal T_{\rho,s}\}
 =
 \ell_\rho(s),
\end{equation}
and
\begin{equation}\label{eq:unique-outer-contact}
 \mathcal T_{\rho,s}\cap\ell_\rho(s)\T
 =
 \{\ell_\rho(s)\}.
\end{equation}
The outer point is generated by the unique Davis--Wielandt pair
$(z_{\rho,s},b_{\rho,s})$ in~\eqref{eq:contact-pair}, which is realized
by $N_\rho$.
\end{theorem}

The proof rests on a scalar support estimate for the exact region
$\mathfrak D_\rho$. Besides controlling the phase in the deformed
numerical-range map, the lemma identifies the unique
Davis--Wielandt pair that can produce the outer contact.

\begin{lemma}\label{lem:phase-support}
Let $1<\rho\leq2$, $0<s<1$, and $t=\ell_\rho(s)$. For
$(z,b)\in\mathfrak D_\rho$, put
\[
 d=1+2s\Real z+s^2b,\qquad
 n=s(1+b)+z+s^2\overline z,\qquad
 m=s^2+2s\Real z+b.
\]

\noindent{\rm (i)} For every $\zeta\in\T$,
\begin{equation}\label{eq:phase-functional}
 \mathcal F_\zeta(z,b)
 :=
 \rho t^2d-2pt\Real(\overline\zeta n)-qm
 \geq0.
\end{equation}
The inequality is strict whenever $\zeta\neq1$.

\noindent{\rm (ii)} Set
\begin{equation}\label{eq:gamma-kappa}
 \gamma=(\rho s-p)t-qs,
 \qquad
 \kappa=q(1+s^2)+2pst-\rho s^2t.
\end{equation}
Then $\kappa>0$ and
\begin{equation}\label{eq:contact-support}
 \mathcal F_1(z,b)
 =
 (1-s^2)
 \bigl(\rho t+2\gamma\Real z-\kappa b\bigr)
 \geq0.
\end{equation}
Equality holds exactly at
\begin{equation}\label{eq:contact-pair}
 z_{\rho,s}
 =
 \begin{cases}
 -\sgn(\gamma)a_{\rho,s},&\gamma\neq0,\\
 0,&\gamma=0,
 \end{cases}
 \qquad
 b_{\rho,s}=\beta_\rho(a_{\rho,s}),
\end{equation}
where
\begin{equation}\label{eq:a-rho-s}
 a_{\rho,s}
 =
 \frac{\rho\abs{\gamma}}
 {\sqrt{4\gamma^2+\rho^2\kappa^2}}.
\end{equation}
Moreover, $a_{\rho,s}<a_\rho$ if $\rho<2$, while
$a_{2,s}<1$.
\end{lemma}

\begin{proof}
Equation~\eqref{eq:ell-quadratic} gives $t>1$, since its left-hand side
at $t=1$ equals $-2ps(1-s)<0$.

We first prove~{\rm (i)} for $\zeta\neq1$. Write
$c=\Real\zeta$, $\delta=1-s^2$, and
\[
 \mathcal F_\zeta(z,b)=A+2\Real(Gz)+Eb.
\]
A direct expansion, using~\eqref{eq:ell-quadratic}, gives
\begin{align*}
 A&=\rho t\delta+2pst(1-c),\\
 G&=sA-\delta(qs+pt\overline\zeta),\\
 E&=s^2A-\delta\bigl(q(1+s^2)+2pstc\bigr),
\end{align*}
together with
\begin{equation}\label{eq:phase-delta}
 \frac{A^2}{\rho^2}+AE-\abs{G}^2
 =
 \left(\frac A\rho-\delta t\right)
 \left(\frac A\rho+\delta t\right)>0.
\end{equation}
Indeed, since $\zeta\neq1$, one has $c<1$, and therefore
\[
 \frac A\rho-\delta t
 =
 \frac{2pst}{\rho}(1-c)>0.
\]
In particular, $A>0$. Set $u=\abs{G}$.

We shall also need the inequality $A+E>2u$. For $\lambda\in\T$, direct
substitution gives
\[
 A+2\Real(G\lambda)+E
 =
 \abs{1+s\lambda}^2
 \left[
 \rho t^2
 -2pt\Real\left(
 \overline\zeta\frac{s+\lambda}{1+s\lambda}
 \right)-q
 \right].
\]
Since
\[
 \left|\frac{s+\lambda}{1+s\lambda}\right|=1,
\]
the term in brackets is bounded below by
\[
 \rho t^2-2pt-q
 =
 (t-1)(\rho t+q)>0.
\]
Moreover, $\abs{1+s\lambda}\geq1-s$, and hence
\[
 A+2\Real(G\lambda)+E
 \geq
 (1-s)^2(t-1)(\rho t+q)>0.
\]
Taking the minimum over $\lambda\in\T$ gives
\[
 A+E-2u>0.
\]

Let now $a=\abs{z}$. Since
$a^2\leq b\leq\beta_\rho(a)$ and
$\Real(Gz)\geq-u a$, we distinguish the sign of $E$.

Suppose first that $E\geq0$. Then
\[
 \mathcal F_\zeta(z,b)
 \geq
 A-2ua+Ea^2.
\]
If $u\geq E$, the function $2ua-Ea^2$ is increasing on $[0,1]$, so
\[
 2ua-Ea^2\leq2u-E<A,
\]
where the last inequality follows from $A+E>2u$. If $u<E$, its maximum
on $[0,1]$ is $u^2/E$, and therefore
\[
 2ua-Ea^2
 \leq
 \frac{u^2}{E}
 <E<A.
\]
Here $E<A$ follows from
\[
 A-E=\delta(\rho t^2+q)>0.
\]
Thus $\mathcal F_\zeta(z,b)>0$ whenever $E\geq0$.

Suppose now that $E<0$, and put $k=-E>0$. Since
$b\leq\beta_\rho(a)$,
\[
 \mathcal F_\zeta(z,b)
 \geq
 A-\bigl(2ua+k\beta_\rho(a)\bigr).
\]
If $qu\geq pk$, Lemma~\ref{lem:beta-support} gives
\[
 2ua+k\beta_\rho(a)\leq2u+k.
\]
Since $A+E>2u$ and $E=-k$, we have $2u+k<A$, and hence
$\mathcal F_\zeta(z,b)>0$.

It remains to consider $qu\leq pk$. Set
\[
 \Delta
 =
 \frac{A^2}{\rho^2}+AE-u^2>0.
\]
Since $E=-k$, this can be written as
\[
 \Delta
 =
 A\left(\frac{A}{\rho^2}-k\right)-u^2,
\]
so $A/\rho^2-k>0$. Moreover,
\[
 \left(\frac{2A}{\rho}-\rho k\right)^2
 -\bigl(\rho^2k^2+4u^2\bigr)
 =
 4\Delta>0.
\]
Consequently,
\[
 \sqrt{\rho^2k^2+4u^2}
 <
 \frac{2A}{\rho}-\rho k,
\]
and therefore the first value in
Lemma~\ref{lem:beta-support} satisfies
\[
 \frac{\rho}{2}
 \left(
 \rho k+\sqrt{\rho^2k^2+4u^2}
 \right)
 <A.
\]
Hence $\mathcal F_\zeta(z,b)>0$ also in this case. This proves
{\rm (i)}.

We now prove~{\rm (ii)}. At $\zeta=1$, direct expansion gives
\begin{align*}
 \mathcal F_1(z,b)
 &=
 \rho t^2-2pst-qs^2\\
 &\quad
 +2\bigl(\rho st^2-pt(1+s^2)-qs\bigr)\Real z\\
 &\quad
 +\bigl(\rho s^2t^2-2pst-q\bigr)b.
\end{align*}
Using~\eqref{eq:ell-quadratic} in each coefficient yields
\[
 \mathcal F_1(z,b)
 =
 (1-s^2)
 \bigl(\rho t+2\gamma\Real z-\kappa b\bigr),
\]
which is~\eqref{eq:contact-support}.

The same quadratic relation also gives
\begin{equation}\label{eq:gamma-kappa-identities}
 t^2-\rho t\kappa-\gamma^2=0,
 \qquad
 \frac{\gamma}{t}
 =
 p+\frac{\rho(1-t)}{s},
 \qquad
 \gamma+pt=s(\rho t-q).
\end{equation}
Since $t>1$, the second identity gives $\gamma<pt$, while the third
gives $\gamma>-pt$. Hence
\[
 -pt<\gamma<pt.
\]
The first identity in~\eqref{eq:gamma-kappa-identities} now gives
\[
 \kappa
 =
 \frac{t^2-\gamma^2}{\rho t}>0.
\]

Set $v=\abs{\gamma}/t$. Then $v<p$, and a direct calculation gives
\[
 p\kappa-q\abs{\gamma}
 =
 \frac{t}{\rho}(p-v)(1+pv)>0.
\]
Thus
\[
 q\abs{\gamma}<p\kappa,
\]
so the first case of Lemma~\ref{lem:beta-support} applies with
$u=\abs{\gamma}$ and $k=\kappa$. It follows that
\[
 2\abs{\gamma}a+\kappa\beta_\rho(a)
 \leq
 \frac{\rho}{2}
 \left(
 \rho\kappa+
 \sqrt{\rho^2\kappa^2+4\gamma^2}
 \right).
\]
From the first identity in~\eqref{eq:gamma-kappa-identities},
\[
 \sqrt{\rho^2\kappa^2+4\gamma^2}
 =
 2t-\rho\kappa,
\]
and therefore
\begin{equation}\label{eq:contact-max}
 2\abs{\gamma}a+\kappa\beta_\rho(a)
 \leq
 \rho t.
\end{equation}
The maximizer is unique and, by
Lemma~\ref{lem:beta-support}, equals
\[
 a_{\rho,s}
 =
 \frac{\rho\abs{\gamma}}
 {\sqrt{4\gamma^2+\rho^2\kappa^2}}.
\]
Since $q\abs{\gamma}<p\kappa$, it lies strictly on the curved branch:
$a_{\rho,s}<a_\rho$ when $\rho<2$, and $a_{2,s}<1$.

Finally, with $a=\abs{z}$,
\[
 \kappa b-2\gamma\Real z
 \leq
 \kappa\beta_\rho(a)+2\abs{\gamma}a
 \leq
 \rho t.
\]
This proves~\eqref{eq:contact-support}. Equality requires equality in
both inequalities above. By the uniqueness in
\eqref{eq:contact-max}, this forces
$a=a_{\rho,s}$ and $b=\beta_\rho(a_{\rho,s})$, together with
\[
 \gamma\Real z=-\abs{\gamma}a_{\rho,s}.
\]
If $\gamma\neq0$, the latter condition and $\abs{z}=a_{\rho,s}$ give
$z=-\sgn(\gamma)a_{\rho,s}$. If $\gamma=0$, the unique maximizer is
$a_{\rho,s}=0$, hence $z=0$. These are exactly the pairs in
\eqref{eq:contact-pair}.
\end{proof}

We now return to the proof of Theorem~\ref{thm:radial-teardrop}.

\begin{proof}[Proof of Theorem~\ref{thm:radial-teardrop}]
Put
\[
 S=\Theta_{\rho,s}(\mathfrak D_\rho),
 \qquad
 t=\ell_\rho(s).
\]

\smallskip
\noindent\emph{Step 1: radial bound for the generating set.}
Fix $(z,b)\in\mathfrak D_\rho$ and write
\[
 (Z,B)=\Psi_s(z,b)=\left(\frac nd,\frac md\right).
\]
Since $b\geq\abs{z}^2$, we have
$d\geq\abs{1+sz}^2\geq(1-s)^2>0$.

Assume first that $Z\neq0$, and set
\[
 R=\abs{\Phi_\rho(Z,B)},
 \qquad
 \zeta=\frac{Z}{\abs{Z}}.
\]
By~\eqref{eq:Phi-quadratic}, $R$ is the positive root of
\[
 P_{Z,B}(r)
 :=
 \rho r^2-2p\abs{Z}r-qB,
 \qquad r\geq0.
\]
Hence $P_{Z,B}(R)=0$. On the other hand,
\eqref{eq:phase-functional} gives
\[
 dP_{Z,B}(t)
 =
 \mathcal F_\zeta(z,b)
 \geq0,
\]
and therefore $P_{Z,B}(t)\geq0$.

The polynomial $P_{Z,B}$ has exactly one positive zero, namely $R$, and
is negative on $(0,R)$. Consequently, $t\geq R$. If $Z=0$, then
\eqref{eq:Phi-rho} gives $\Phi_\rho(Z,B)=0$. Thus
\begin{equation}\label{eq:S-radial-bound}
 S\subseteq t\overline\D.
\end{equation}

\smallskip
\noindent\emph{Step 2: equality among the generating points.}
Suppose
$\lambda=\Theta_{\rho,s}(z,b)\in S$ and $\abs{\lambda}=t$.
Then $Z\neq0$ and $R=t$. Hence
\[
 \mathcal F_\zeta(z,b)
 =
 dP_{Z,B}(t)
 =
 0.
\]
Lemma~\ref{lem:phase-support} therefore forces
\[
 \zeta=1,
 \qquad
 (z,b)=(z_{\rho,s},b_{\rho,s}).
\]
Thus a generating point can lie on $t\T$ only at the unique
Davis--Wielandt pair in~\eqref{eq:contact-pair}.

Conversely, evaluate at this pair. By
Lemma~\ref{lem:phase-support},
$\mathcal F_1=0$, whereas $\mathcal F_\xi>0$ for every
$\xi\neq1$. Since
\[
 \mathcal F_\xi-\mathcal F_1
 =
 2pt\bigl(
 \Real n-\Real(\overline\xi n)
 \bigr),
\]
the function
$\xi\mapsto\Real(\overline\xi n)$ on $\T$ has its unique maximum at
$\xi=1$. Hence $n$ is a positive real number. Since $d>0$, it follows
that $Z=n/d>0$, and therefore
\[
 P_{Z,B}(t)=0
 \quad\Longrightarrow\quad
 \Theta_{\rho,s}(z_{\rho,s},b_{\rho,s})=t.
\]
Consequently,
\begin{equation}\label{eq:S-outer-contact}
 S\cap t\T=\{t\}.
\end{equation}
The pair $(z_{\rho,s},b_{\rho,s})$ lies on the curved upper boundary of
$\mathfrak D_\rho$ by Lemma~\ref{lem:phase-support}, and the construction
in Theorem~\ref{thm:exact-DW} realizes it by $N_\rho$.

\smallskip
\noindent\emph{Step 3: passage to the closure of the generating set.}
Let
\[
 \lambda_j
 =
 \Theta_{\rho,s}(z_j,b_j)\in S,
 \qquad
 \lambda_j\longrightarrow\lambda\in t\T.
\]
Set $R_j=\abs{\lambda_j}$. Then $R_j\to t$, and the corresponding
$Z_j$ is nonzero for all sufficiently large $j$. Put
$\zeta_j=Z_j/\abs{Z_j}$. Since
$P_{Z_j,B_j}(R_j)=0$, subtraction from the phase inequality at $t$
gives
\begin{equation}\label{eq:phase-limit-identity}
 \mathcal F_{\zeta_j}(z_j,b_j)
 =
 d_j(t-R_j)
 \bigl(
 \rho(t+R_j)-2p\abs{Z_j}
 \bigr).
\end{equation}
The set $\mathfrak D_\rho$ is compact. Moreover,
$d_j\geq(1-s)^2$, while $d_j$, $Z_j=n_j/d_j$, and $B_j=m_j/d_j$
remain bounded. Hence the factor multiplying $t-R_j$ on the
right-hand side of~\eqref{eq:phase-limit-identity} is bounded, and
therefore
\[
 \mathcal F_{\zeta_j}(z_j,b_j)\longrightarrow0.
\]

After passing to a subsequence, compactness of
$\mathfrak D_\rho\times\T$ gives
\[
 (z_j,b_j,\zeta_j)
 \longrightarrow
 (z,b,\zeta).
\]
By continuity,
$\mathcal F_\zeta(z,b)=0$. Lemma~\ref{lem:phase-support} then yields
\[
 \zeta=1,
 \qquad
 (z,b)=(z_{\rho,s},b_{\rho,s}).
\]
At this pair we have $n>0$, so $\Theta_{\rho,s}$ is continuous in a
neighborhood of $(z_{\rho,s},b_{\rho,s})$. The chosen subsequence
therefore converges to
\[
 \Theta_{\rho,s}(z_{\rho,s},b_{\rho,s})=t.
\]
Since the original sequence converges to $\lambda$, we obtain
$\lambda=t$. Hence
\begin{equation}\label{eq:closure-outer-contact}
 \overline S\cap t\T=\{t\}.
\end{equation}

\smallskip
\noindent\emph{Step 4: passage to the closed convex hull.}
Let
\[
 \lambda\in\overline{\conv}S\cap t\T,
 \qquad
 \eta=\frac{\lambda}{t}\in\T.
\]
Choose finite convex combinations
\[
 y_k
 =
 \sum_{j=1}^{n_k}\alpha_{k,j}\lambda_{k,j}
 \longrightarrow\lambda,
 \qquad
 \lambda_{k,j}\in S,
\]
where $\alpha_{k,j}\geq0$ and
$\sum_j\alpha_{k,j}=1$. For each $k,j$, define
\[
 \delta_{k,j}
 =
 t-\Real(\overline\eta\lambda_{k,j})\geq0.
\]
Then
\[
 \sum_j\alpha_{k,j}\delta_{k,j}
 =
 t-\Real(\overline\eta y_k)
 \longrightarrow0.
\]
Choose an index $j(k)$ with $\alpha_{k,j(k)}>0$ such that
$\delta_{k,j(k)}$ does not exceed this weighted average. Hence
$\delta_{k,j(k)}\to0$.

By~\eqref{eq:S-radial-bound},
$\lambda_{k,j(k)}\in t\overline\D$. Passing to a subsequence, let
\[
 \lambda_{k,j(k)}\longrightarrow\mu\in t\overline\D.
\]
The preceding limit gives
\[
 \Real(\overline\eta\mu)=t.
\]
Since
\[
 \Real(\overline\eta\mu)
 \leq
 \abs{\mu}
 \leq
 t,
\]
equality must hold throughout, and therefore
$\mu=t\eta=\lambda$. Thus $\lambda\in\overline S\cap t\T$, and
\eqref{eq:closure-outer-contact} yields $\lambda=t$.

Finally, Step~2 gives $t\in S$, while Step~1 gives
\[
 \mathcal T_{\rho,s}
 =
 \overline{\conv}S
 \subseteq
 t\overline\D.
\]
Hence the outer radius of $\mathcal T_{\rho,s}$ is $t=\ell_\rho(s)$,
and its intersection with $t\T$ consists only of $t$. This proves
\eqref{eq:teardrop-radius} and~\eqref{eq:unique-outer-contact}.
\end{proof}

When $s=0$, the universal set reduces to $\overline\D$. Indeed,
\eqref{eq:deformed-facts} gives one inclusion, while scalar contractions
give the reverse inclusion. In this case the outer contact is not
unique.

\subsection{The ordinary numerical-range companion}

The ordinary projection has a simpler form. Write
$\overline D(a,r)=\{z:|z-a|\leq r\}$.

\begin{theorem}\label{thm:ordinary-teardrop}
For $1\leq\rho\leq2$ and $0\leq s<1$,
\begin{equation}\label{eq:ordinary-teardrop}
 \overline{\bigcup_{C\in\mathcal C_\rho}W(g_s(C))}
 =\conv\left(\overline\D\cup
 \overline D\left(s,\frac\rho2(1-s^2)\right)\right).
\end{equation}
\end{theorem}

\begin{proof}
The case $\rho=1$ follows from von Neumann's inequality
(see~\cite{SzNagyFoias}), and the second disk is contained in
$\overline{\D}$. Assume $1<\rho\leq2$, put
$X=g_s(C)$ and $D=I+sC$, and fix $e^{i\theta}\in\T$. The support function
of the set on the right of~\eqref{eq:ordinary-teardrop} is
\begin{equation}\label{eq:ordinary-support}
 h(\theta)=\max\left\{1,
 s\cos\theta+\frac\rho2(1-s^2)\right\}.
\end{equation}

Suppose first that the first term is maximal. Set
\[
 \eta=\frac{e^{i\theta}-s}{1-se^{i\theta}},\qquad
 \delta=|e^{i\theta}-s|^2.
\]
Then $|\eta|=1$, and direct expansion gives
\begin{align}\label{eq:unit-support-factorization}
 D^*\bigl(I-\Real(e^{-i\theta}X)\bigr)D
 &=\frac{1-s^2}{2}Q_C(\eta,1)\\
 &\quad+\frac{\delta-p(1-s^2)}{2}
 (C-\eta I)^*(C-\eta I).\notag
\end{align}
The second coefficient is nonnegative exactly when
$q+\rho s^2-2s\cos\theta\geq0$, which is equivalent to the first regime
in~\eqref{eq:ordinary-support}. Thus the left side is positive.

In the second regime put
\[
 u_\theta=\frac{e^{i\theta}-\rho s}{p}.
\]
The regime condition is equivalent to $|u_\theta|\leq1$. A second direct
expansion gives
\begin{align}\label{eq:disk-support-factorization}
 D^*&\left[\left(s\cos\theta+\frac\rho2(1-s^2)\right)I
 -\Real(e^{-i\theta}X)\right]D\\
 &=\frac{\rho(1-s^2)}2P_C(u_\theta)\geq0.\notag
\end{align}
Since $D$ is invertible,~\eqref{eq:unit-support-factorization} and
\eqref{eq:disk-support-factorization} prove the inclusion from left to
right in every half-plane.

For the reverse inclusion, let $\lambda_1,\ldots,\lambda_n\in\T$ and
let $\alpha_1,\ldots,\alpha_n\geq0$ satisfy
$\sum_j\alpha_j=1$. Since $g_s$ maps $\T$ onto itself, choose
$\eta_j\in\T$ such that $g_s(\eta_j)=\lambda_j$. The scalar operators
$\eta_j I$ belong to $\mathcal C_\rho$, and the direct sum
$C=\bigoplus_{j=1}^n\eta_j I$ satisfies
\[
 \sum_{j=1}^n\alpha_j\lambda_j
 =
 \inner{g_s(C)x}{x},
 \qquad
 x=(\sqrt{\alpha_1},\ldots,\sqrt{\alpha_n})^{\mathsf T}.
\]
Thus the closed union on the left of~\eqref{eq:ordinary-teardrop}
contains $\overline\D$. Moreover,
\[
 g_s(N_\rho)=
 \begin{pmatrix}s&\rho(1-s^2)\\0&s\end{pmatrix},
\qquad
 W(g_s(N_\rho))=\overline D\left(s,\frac\rho2(1-s^2)\right).
\]
Taking direct sums of $N_\rho$ with the preceding scalar operators shows
that the closed union contains every finite convex combination of points
from the two disks, and hence it contains their convex hull. This proves
the reverse inclusion. At $\rho=2$ the formula is Drury's teardrop,
while at $\rho=1$ it reduces to the unit disk.
\end{proof}

\begin{example}
Let $\rho=3/2$ and $s=1/2$. The ordinary companion is
\[
\conv\left(
\overline{\D}\cup
\overline D\left(\frac12,\frac{9}{16}\right)
\right),
\]
whose rightmost point is $17/16$. In contrast,
\[
\max_{\lambda\in\mathcal T_{3/2,\,1/2}}\abs{\lambda}
=
\ell_{3/2}\left(\frac12\right)
=
\frac{13+\sqrt{217}}{24}
\approx1.1554,
\]
and the outer circle meets $\mathcal T_{3/2,\,1/2}$ only at this
positive real point.
\end{example}

\begin{figure}[htbp]
\centering
\input{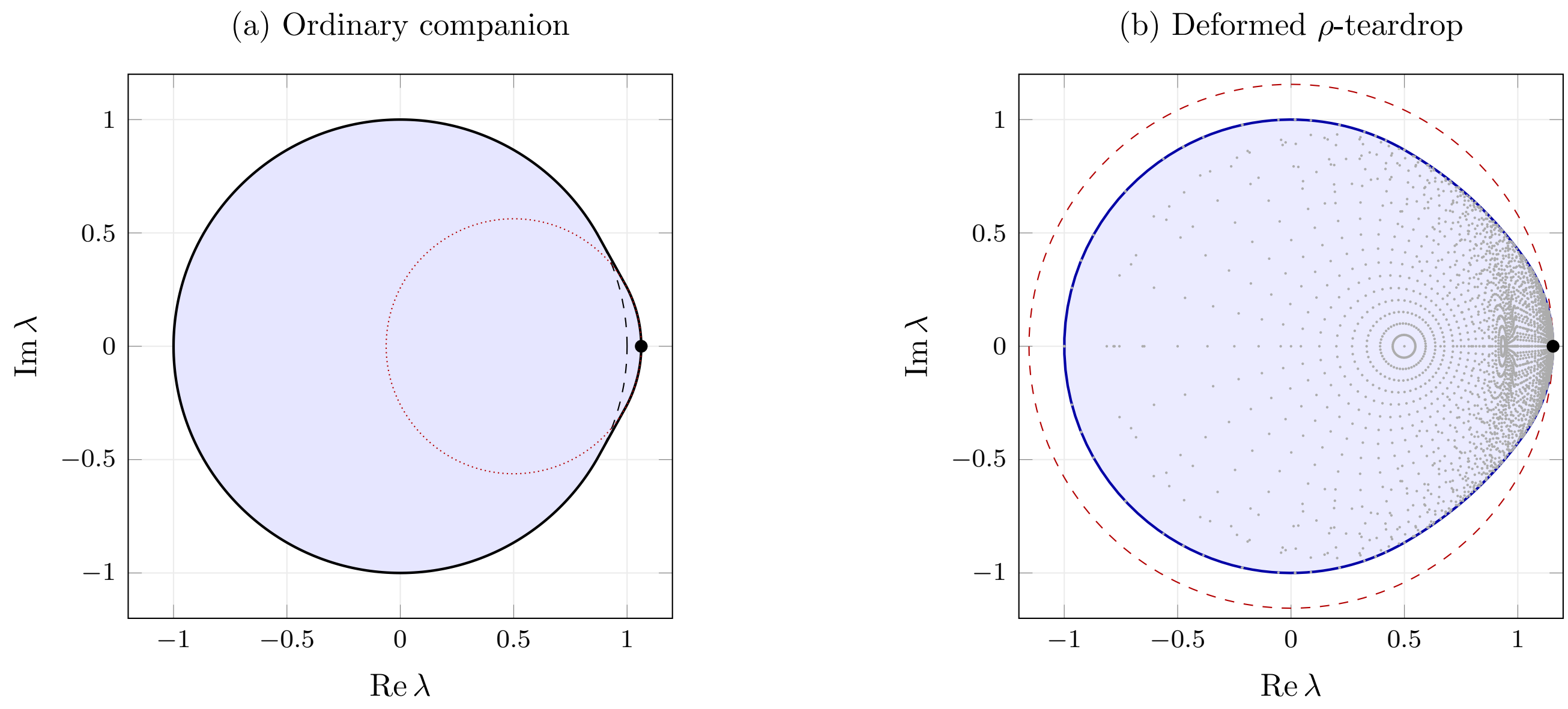}

\begin{minipage}{.485\textwidth}
\centering
\begin{tikzpicture}
\begin{axis}[
 width=\textwidth,height=.82\textwidth,
 axis equal image,
 xmin=-1.20,xmax=1.20,
 ymin=-1.20,ymax=1.20,
 xlabel={$\Real\lambda$},
 ylabel={$\operatorname{Im}\lambda$},
 title={(a) Ordinary companion},
 title style={font=\small},
 label style={font=\small},
 tick label style={font=\scriptsize},
 grid=major,
 grid style={gray!15}
]
\addplot[draw=black,thick,fill=blue!10]
 table{figures/ordinary_hull_rho1p5_s0p5.dat};

\addplot[black,dashed,domain=0:360,samples=181]
 ({cos(x)},{sin(x)});

\addplot[red!70!black,densely dotted,domain=0:360,samples=181]
 ({0.5+0.5625*cos(x)},{0.5625*sin(x)});

\addplot[only marks,mark=*,mark size=1.8pt]
 coordinates {(1.0625,0)};
\end{axis}
\end{tikzpicture}
\end{minipage}\hfill
\begin{minipage}{.485\textwidth}
\centering
\begin{tikzpicture}
\begin{axis}[
 width=\textwidth,height=.82\textwidth,
 axis equal image,
 xmin=-1.20,xmax=1.20,
 ymin=-1.20,ymax=1.20,
 xlabel={$\Real\lambda$},
 ylabel={$\operatorname{Im}\lambda$},
 title={(b) Deformed $\rho$-teardrop},
 title style={font=\small},
 label style={font=\small},
 tick label style={font=\scriptsize},
 grid=major,
 grid style={gray!15}
]
\addplot[draw=blue!65!black,thick,fill=blue!8]
 table{figures/deformed_hull_rho1p5_s0p5.dat};

\addplot[only marks,mark=*,mark size=.20pt,gray!65]
 table{figures/deformed_generators_rho1p5_s0p5.dat};

\addplot[red!70!black,dashed,domain=0:360,samples=181]
 ({\figradius*cos(x)},{\figradius*sin(x)});

\addplot[only marks,mark=*,mark size=1.8pt]
 coordinates {(\figradius,0)};
\end{axis}
\end{tikzpicture}
\end{minipage}

\caption{Graphical illustration of the preceding example for
$\rho=3/2$ and $s=1/2$. In {\rm(a)}, the dashed unit circle and the
dotted circle
$\partial D(s,\rho(1-s^2)/2)$ generate the ordinary companion.
In {\rm(b)}, the gray points sample the deformed generators and the
solid boundary is their convex hull $\mathcal T_{\rho,s}$. The dashed
circle has radius $\ell_\rho(s)$ and meets $\mathcal T_{\rho,s}$ only
at the positive real point $\ell_\rho(s)$.}
\label{fig:teardrop-comparison}
\end{figure}

\FloatBarrier
\section{The sharp functional-calculus theorem}

The universal geometry developed above yields the sharp M\"obius
estimate directly.

\begin{corollary}\label{cor:Mobius-bound}
Let $1\leq\rho\leq2$, $C\in\mathcal C_\rho$, and $0\leq s<1$. Then
\[
 w_\rho(g_s(C))\leq\ell_\rho(s).
\]
For $1<\rho\leq2$, equality is attained by $C=N_\rho$.
\end{corollary}

\begin{proof}
For $1<\rho\leq2$ and $s>0$, combine
Theorems~\ref{thm:universal-teardrop} and
\ref{thm:radial-teardrop} with~\eqref{eq:gauge-Wrho}.
When $s=0$, one has $g_0(C)=C$, so the estimate follows from
$C\in\mathcal C_\rho$.

For equality at $s=0$,~\eqref{eq:N-pencil} gives
$w_\rho(N_\rho)\leq1$, while $\norm{N_\rho}=\rho$ and the dilation
relation imply
\[
 w_\rho(N_\rho)\geq\frac{\norm{N_\rho}}{\rho}=1.
\]
Thus $w_\rho(N_\rho)=1=\ell_\rho(0)$.
For $1<\rho\leq2$ and $s>0$, Theorem~\ref{thm:radial-teardrop}
shows that the outer contact is realized by $N_\rho$, giving equality.

When $\rho=1$, the estimate follows from the classical von Neumann
inequality; see~\cite{SzNagyFoias}.
\end{proof}

We can now return to the original functional-calculus problem.

\begin{theorem}\label{thm:main}
Let $1\leq\rho\leq2$, let $T\in\BH$ satisfy
\[
 w_\rho(T)\leq1,
\]
and let $f\in A(\overline\D)$ satisfy
\[
 \norm{f}_{\infty,\D}\leq1.
\]
Then
\begin{equation}\label{eq:main-theorem}
 w_\rho(f(T))
 \leq
 \ell_\rho\bigl(\abs{f(0)}\bigr).
\end{equation}
The bound is optimal for every prescribed value of $\abs{f(0)}$.
\end{theorem}

\begin{proof}
Set $s=\abs{f(0)}$. Multiplying $f$ by a unimodular constant does not
change $w_\rho(f(T))$, so we may assume that $f(0)=s$.

If $s=1$, the maximum modulus principle implies that $f$ is constant,
and~\eqref{eq:main-theorem} follows immediately. Assume therefore that
$s<1$ and define
\[
 h(z)=\frac{f(z)-s}{1-sf(z)}.
\]
For $z\in\overline\D$,
\[
 \abs{1-sf(z)}
 \geq
 1-s\abs{f(z)}
 \geq
 1-s
 >
 0.
\]
Thus $h\in A(\overline\D)$. Moreover, $h(0)=0$, and the disk
automorphism $z\mapsto(z-s)/(1-sz)$ gives
$\norm{h}_{\infty,\D}\leq1$.

Since $w_\rho(T)\leq1$, Lemma~\ref{lem:positivity}{\rm(i)} gives
$T\in\mathcal C_\rho$, and Theorem~\ref{thm:zero-mapping} yields
\[
 C:=h(T)\in\mathcal C_\rho.
\]
Since the disk-algebra functional calculus is a unital homomorphism and
$f=g_s\circ h$, one has
\[
 f(T)=g_s(C).
\]
Corollary~\ref{cor:Mobius-bound} therefore proves
\eqref{eq:main-theorem}.

It remains to prove sharpness. Since $\norm{N_\rho}=\rho$, the defining
dilation relation implies
$w_\rho(N_\rho)\geq\norm{N_\rho}/\rho=1$, while
\eqref{eq:N-pencil} gives $N_\rho\in\mathcal C_\rho$ and hence
$w_\rho(N_\rho)\leq1$. Thus
\[
 w_\rho(N_\rho)=1.
\]

For $1<\rho\leq2$ and $0<s<1$, take $T=N_\rho$ and $f=g_s$. The
function $g_s$ belongs to $A(\overline\D)$.
Theorem~\ref{thm:radial-teardrop} gives
\[
 w_\rho(g_s(N_\rho))=\ell_\rho(s),
\]
so equality holds. When $s=0$, take again $T=N_\rho$ and $f=g_0$.
Since $g_0(z)=z$, equality follows from
$w_\rho(N_\rho)=1=\ell_\rho(0)$. For $s=1$, take the constant function
$f=1$.

Finally, when $\rho=1$ and $0\leq s<1$, take the scalar operator
$T=I$ and $f=g_s$. Since $g_s(1)=1$,
\[
 w_1(f(T))=1=\ell_1(s).
\]
For $s=1$, again take $f=1$. Thus the bound is optimal for every
prescribed value of $\abs{f(0)}$.
\end{proof}

For $0\leq s\leq1$, define the universal optimal pointwise constant by
\[
M_\rho(s)
:=
\sup\left\{
w_\rho(f(T)) :
T\in\mathcal B(\mathcal H),\;
w_\rho(T)\leq1,\;
f\in A(\overline{\D}),\;
\norm{f}_{\infty,\D}\leq1,\;
\abs{f(0)}=s
\right\},
\]
Theorem~\ref{thm:main} and its extremals give the exact formula
\begin{equation}\label{eq:M-rho}
 M_\rho(s)=\ell_\rho(s),
 \qquad 0\leq s\leq1.
\end{equation}
Maximizing this expression recovers the sharp uniform constant
in~\eqref{eq:uniform-intro}.

\begin{proposition}\label{prop:uniform-maximum}
For $1<\rho\leq2$, the function $\ell_\rho$ has a unique maximum at
\begin{equation}\label{eq:s-star}
 s_\rho^*
 =
 \frac{\sqrt{5\rho^2+2\rho+1}-(\rho+1)}{2\rho},
\end{equation}
and
\[
 \ell_\rho(s_\rho^*)=C_\rho.
\]
For $\rho=1$, one has $\ell_1\equiv1$.
\end{proposition}

\begin{proof}
Write $t=\ell_\rho(s)$. Differentiating
\eqref{eq:ell-quadratic},
\[
 \rho t^2-B_\rho(s)t-qs^2=0,
\]
gives
\[
 \bigl(2\rho t-B_\rho(s)\bigr)t'
 =
 B_\rho'(s)t+2qs.
\]
Since
\[
 B_\rho'(s)=2(p-\rho s)
\]
and
\[
 2\rho t-B_\rho(s)
 =
 \sqrt{B_\rho(s)^2+4\rho qs^2}>0,
\]
an interior critical point satisfies
\[
 (\rho s-p)t=qs.
\]
Substituting this relation into~\eqref{eq:ell-quadratic} reduces the
critical-point equation to
\[
 \rho s^2+(\rho+1)s-\rho=0.
\]
Its unique root in $(0,1)$ is precisely $s_\rho^*$ in
\eqref{eq:s-star}.

At $s=0$, one has $t=1$, and the differentiated identity gives
$\ell_\rho'(0)>0$. At $s=1$, again $t=1$, and it gives
$\ell_\rho'(1)<0$. Since there is only one critical point in $(0,1)$,
$s_\rho^*$ is the unique maximizer.

Combining
\[
 (\rho s-p)t=qs
\]
with
\[
 \rho s^2+(\rho+1)s-\rho=0
\]
gives, at the critical point,
\[
 t=1+\frac{ps}{\rho}.
\]
Hence
\[
 \ell_\rho(s_\rho^*)
 =
 1+\frac{p}{\rho}s_\rho^*
 =
 C_\rho.
\]
Finally, $\ell_1(s)=1$ follows directly from~\eqref{eq:B-ell}.
\end{proof}

At the endpoints,
\[
 \ell_1(s)=1,
 \qquad
 \ell_2(s)=1+s-s^2.
\]
Thus Theorem~\ref{thm:main} interpolates between von Neumann's
inequality~\cite{SzNagyFoias} and Drury's estimate~\cite{Drury}.

\section{Extremal contact and equality}

Extremal numerical-radius and $\rho$-contraction problems often exhibit
rigidity through the appearance of canonical finite-dimensional models
on reducing subspaces; see
\cite{HaagerupDeLaHarpe,GauWang,GauWu2026}.
In our setting, this rigidity is produced by the curved boundary of the
universal Davis--Wielandt region.

We begin with the corresponding boundary rigidity result.

\begin{lemma}\label{lem:curved-rigidity}
Let $1<\rho\leq2$, let $C\in\mathcal C_\rho$, and let $x$ be a unit
vector. Put
$a=\abs{\inner{Cx}{x}}$ and $b=\norm{Cx}^2$. Assume
\begin{equation}\label{eq:curved-equality}
 b
 =
 \frac{\rho}{2}
 \left(
 \rho+\sqrt{\rho^2-4a^2}
 \right),
\end{equation}
where $a<a_\rho$ if $\rho<2$, and $a<1$ if $\rho=2$. Then $C$ has a
reducing two-dimensional subspace on which it is unitarily equivalent
to $N_\rho$.
\end{lemma}

\begin{proof}
Multiplying $C$ by a unimodular scalar preserves membership in
$\mathcal C_\rho$ and replaces $N_\rho$ only by a unitarily equivalent
matrix. We may therefore assume that
$\inner{Cx}{x}=a$.

Set
\[
 r_0=\frac{\rho a}{pb}.
\]
The strict curved-branch assumption gives $0\leq r_0<1$. Moreover,
\eqref{eq:curved-equality} is equivalent to
\[
 \frac{b}{\rho^2}+\frac{a^2}{b}=1.
\]
Substituting $r_0$ into the quadratic form of $P_C(r_0)$ therefore gives
$\inner{P_C(r_0)x}{x}=0$. By
Lemma~\ref{lem:positivity}{\rm(ii)}, $P_C(r_0)\geq0$, and a positive
operator with vanishing quadratic form at $x$ must annihilate $x$.
Hence
\begin{equation}\label{eq:P-kernel}
 P_C(r_0)x=0.
\end{equation}

Let $S=S_C$ and $F=F_C$ be as in
\eqref{eq:positive-real-factorization}, and set
$v=S(r_0)^{-1}x$. Since $S(r_0)v=x$, the factorization and
\eqref{eq:P-kernel} give
\[
 \inner{\Real F(r_0)v}{v}=0.
\]
The scalar function
\[
 w\longmapsto\inner{F(w)v}{v}
\]
is analytic on $\D$ and has nonnegative real part. Its real part is
therefore a nonnegative harmonic function attaining its minimum zero at
the interior point $r_0$. By the minimum principle, it vanishes
identically on $\D$.

Since $\Real F(w)\geq0$, we obtain
$\inner{\Real F(w)v}{v}=0$ for every $w\in\D$, and hence
\[
 \Real F(w)v=0.
\]
Equivalently,
\[
 F(w)v+F(w)^*v=0.
\]
The first term is holomorphic in $w$, whereas the second is
anti-holomorphic. Applying $\partial/\partial w$ to this identity gives
\[
F'(w)v=0.
\] Since
\[
 F'(w)
 =
 \frac{2p}{\rho^2}C^2S(w)^2,
\]
evaluation at $w=r_0$ yields
\[
 C^2S(r_0)x=0.
\]
The operators $C$ and $S(r_0)$ commute, and $S(r_0)$ is invertible.
Consequently,
\[
 C^2x=0.
\]

We now identify the resulting invariant subspace. Since the strict
curved-boundary condition gives $b>a^2$, the vectors $x$ and $Cx$ are
linearly independent. Thus
\[
 \mathcal M=\operatorname{span}\{x,Cx\}
\]
is two-dimensional, and $C^2x=0$ shows that it is invariant under $C$.
Put
\[
 d=\sqrt{b-a^2},
 \qquad
 y=\frac{Cx-ax}{d}.
\]
Then $(x,y)$ is an orthonormal basis of $\mathcal M$, with
\[
 Cx=ax+dy,
 \qquad
 Cy=-\frac{a^2}{d}x-ay.
\]
Hence
\[
 A:=C|_{\mathcal M}
 =
 \begin{pmatrix}
 a&-a^2/d\\
 d&-a
 \end{pmatrix}.
\]
In particular, $A^2=0$ and $A$ has rank one. Moreover,
\[
 \norm{A}^2
 =
 \frac{b^2}{b-a^2}
 =
 \rho^2,
\]
where the last equality follows from
$b/\rho^2+a^2/b=1$. Thus $A$ is a rank-one nilpotent of norm $\rho$,
and is therefore unitarily equivalent to
\[
 N_\rho=
 \begin{pmatrix}
 0&\rho\\
 0&0
 \end{pmatrix}.
\]

It remains to show that $\mathcal M$ is reducing, not merely invariant.
After a unitary change of basis on $\mathcal M$, write
\[
 C=
 \begin{pmatrix}
 N_\rho&B\\
 0&D_0
 \end{pmatrix}
 \quad\text{on}\quad
 \mathcal M\oplus\mathcal M^\perp.
\]
For $\zeta\in\T$, the upper-left block of the positive operator
$Q_C(\zeta,1)$ is $Q_{N_\rho}(\zeta,1)$. By
\eqref{eq:N-pencil},
\[
 \ker Q_{N_\rho}(\zeta,1)
 =
 \operatorname{span}\{k_\zeta\},
 \qquad
 k_\zeta=(p\overline\zeta,1)^{\mathsf T}.
\]
The lower-left block of $Q_C(\zeta,1)$ is
\[
 -B^*(p\zeta I+qN_\rho).
\]
Positivity of $Q_C(\zeta,1)$ forces this block to annihilate the kernel
of its upper-left diagonal block; hence
\[
 B^*(p\zeta I+qN_\rho)k_\zeta=0.
\]
Using $p^2+\rho q=1$,
\[
 (p\zeta I+qN_\rho)k_\zeta
 =
 \begin{pmatrix}
 1\\
 p\zeta
 \end{pmatrix}.
\]
As $\zeta$ varies over $\T$, these vectors span $\C^2$. Therefore
$B^*=0$, and hence $B=0$. Thus $\mathcal M$ reduces $C$, and the
restriction of $C$ to $\mathcal M$ is unitarily equivalent to
$N_\rho$.
\end{proof}

The supporting hyperplane in Lemma~\ref{lem:phase-support} admits a
natural operator form. For $0<s<1$ and $t=\ell_\rho(s)$, define
\[
 H_{\rho,s}(C)
 =
 \rho tI+\gamma(C+C^*)-\kappa C^*C,
\]
where $\gamma$ and $\kappa$ are given by~\eqref{eq:gamma-kappa}.
The next theorem gives the equality criterion associated with the
unique outer contact.

\begin{theorem}\label{thm:defect-equality}
Let $1<\rho\leq2$, $0<s<1$, $C\in\mathcal C_\rho$, and
$X=g_s(C)$. Then $H_{\rho,s}(C)\geq0$, and the following assertions are
equivalent:
\begin{enumerate}
\item[\rm(i)] $w_\rho(X)=\ell_\rho(s)$;

\item[\rm(ii)] $0\in\sigma(H_{\rho,s}(C))$;

\item[\rm(iii)] there exists a sequence of unit vectors $(x_n)$ such that
\[
 \inner{Cx_n}{x_n}\longrightarrow z_{\rho,s},
 \qquad
 \norm{Cx_n}^2\longrightarrow b_{\rho,s}.
\]
\end{enumerate}
Moreover,
\begin{equation}\label{eq:kernel-reducing}
 \ker H_{\rho,s}(C)\neq\{0\}
 \quad\Longleftrightarrow\quad
 C\simeq N_\rho\oplus C_0
\end{equation}
for some $C_0\in\mathcal C_\rho$.
\end{theorem}

\begin{proof}
Let
\[
 t=\ell_\rho(s),
 \qquad
 D=I+sC.
\]
Since $C\in\mathcal C_\rho$, its spectral radius is at most $1$.
As $0<s<1$, we have $-1/s\notin\sigma(C)$, and hence $D$ is
invertible.

For a unit vector $x$, put
\[
 z=\inner{Cx}{x},
 \qquad
 b=\norm{Cx}^2.
\]
By Theorem~\ref{thm:exact-DW},
$(z,b)\in\mathfrak D_\rho$. Therefore
\eqref{eq:contact-support} gives
\begin{equation}\label{eq:H-quadratic-form}
 \inner{H_{\rho,s}(C)x}{x}
 =
 \rho t+2\gamma\Real z-\kappa b
 \geq0.
\end{equation}
Thus
\[
 H_{\rho,s}(C)\geq0.
\]

\smallskip
\noindent\emph{Equivalence of {\rm(ii)} and {\rm(iii)}.}
Define
\[
 \varphi(z,b)
 =
 \rho t+2\gamma\Real z-\kappa b.
\]
By Lemma~\ref{lem:phase-support}, $\varphi$ is nonnegative on
$\mathfrak D_\rho$, and its zero set consists of the single point
\[
 (z_{\rho,s},b_{\rho,s}).
\]

Suppose first that $(x_n)$ is a sequence of unit vectors satisfying
\[
 \inner{H_{\rho,s}(C)x_n}{x_n}\longrightarrow0.
\]
Set
\[
 z_n=\inner{Cx_n}{x_n},
 \qquad
 b_n=\norm{Cx_n}^2.
\]
Then $(z_n,b_n)\in\mathfrak D_\rho$ and, by
\eqref{eq:H-quadratic-form},
\[
 \varphi(z_n,b_n)\longrightarrow0.
\]
We claim that
\[
 (z_n,b_n)\longrightarrow(z_{\rho,s},b_{\rho,s}).
\]
Indeed, otherwise some subsequence would remain a positive distance
from $(z_{\rho,s},b_{\rho,s})$. Since $\mathfrak D_\rho$ is compact,
that subsequence would admit a further subsequence converging to some
$(z_0,b_0)\in\mathfrak D_\rho$. Continuity of $\varphi$ would then give
\[
 \varphi(z_0,b_0)=0,
\]
while $(z_0,b_0)\neq(z_{\rho,s},b_{\rho,s})$, contradicting the
uniqueness of the zero.

Conversely, if
\[
 \inner{Cx_n}{x_n}\longrightarrow z_{\rho,s},
 \qquad
 \norm{Cx_n}^2\longrightarrow b_{\rho,s},
\]
then continuity of $\varphi$ and
\eqref{eq:H-quadratic-form} give
\[
 \inner{H_{\rho,s}(C)x_n}{x_n}\longrightarrow0.
\]

Finally, since $H_{\rho,s}(C)$ is positive,
\[
 0\in\sigma(H_{\rho,s}(C))
 \quad\Longleftrightarrow\quad
 \inf_{\norm{x}=1}
 \inner{H_{\rho,s}(C)x}{x}=0.
\]
This proves the equivalence of {\rm(ii)} and {\rm(iii)}.

\smallskip
\noindent\emph{{\rm(i)}$\Rightarrow${\rm(ii)}.}
Recall that
\[
 X=g_s(C)=(sI+C)(I+sC)^{-1},
\]
so that
\[
 XD=sI+C.
\]
For every unit vector $y$, put
\[
 z=\inner{Cy}{y},
 \qquad
 b=\norm{Cy}^2.
\]
Using $XD=sI+C$ in the definition of the boundary pencil gives, after
expansion,
\begin{equation}\label{eq:phase-operator-form}
 \inner{D^*Q_X(\zeta,t)Dy}{y}
 =
 \mathcal F_\zeta(z,b),
 \qquad \zeta\in\T.
\end{equation}

For $\zeta=1$, Lemma~\ref{lem:phase-support} gives
\[
 \mathcal F_1(z,b)
 =
 (1-s^2)
 \bigl(\rho t+2\gamma\Real z-\kappa b\bigr).
\]
Comparing with~\eqref{eq:H-quadratic-form}, and since this holds for
every $y$, we obtain the operator identity
\begin{equation}\label{eq:H-congruence}
 D^*Q_X(1,t)D
 =
 (1-s^2)H_{\rho,s}(C).
\end{equation}

Assume {\rm(i)}, so that
\[
 w_\rho(X)=t,
\]
and suppose, to the contrary, that
$0\notin\sigma(H_{\rho,s}(C))$. Since
$H_{\rho,s}(C)\geq0$, there exists $\varepsilon>0$ such that
\[
 H_{\rho,s}(C)\geq\varepsilon I.
\]
By~\eqref{eq:H-congruence},
\[
 D^*Q_X(1,t)D
 \geq
 (1-s^2)\varepsilon I.
\]
Let $u\in\mathcal H$ and put $y=D^{-1}u$. Then
\[
 \begin{aligned}
 \inner{Q_X(1,t)u}{u}
 &=
 \inner{D^*Q_X(1,t)Dy}{y}\\
 &\geq
 (1-s^2)\varepsilon\norm{y}^2\\
 &\geq
 \frac{(1-s^2)\varepsilon}{\norm{D}^2}\norm{u}^2.
 \end{aligned}
\]
Hence $Q_X(1,t)$ is bounded below by a positive multiple of the
identity.

Since the map
\[
 \zeta\longmapsto Q_X(\zeta,t)
\]
is norm-continuous on $\T$, there exist an open arc $U\subset\T$
containing $1$ and a constant $c_1>0$ such that
\[
 Q_X(\zeta,t)\geq c_1I,
 \qquad
 \zeta\in U.
\]

We now control the remaining phases. By part~{\rm(i)} of
Lemma~\ref{lem:phase-support},
\[
 \mathcal F_\zeta(z,b)>0
\]
whenever $\zeta\neq1$ and $(z,b)\in\mathfrak D_\rho$. Therefore, on the
compact set
$
 (\T\setminus U)\times\mathfrak D_\rho,
$
the function $\mathcal F_\zeta(z,b)$ has a strictly positive minimum;
write
\[
 m
 =
 \min_{\substack{\zeta\in\T\setminus U\\
                  (z,b)\in\mathfrak D_\rho}}
 \mathcal F_\zeta(z,b)
 >0.
\]
Equation~\eqref{eq:phase-operator-form} then gives
\[
 D^*Q_X(\zeta,t)D\geq mI,
 \qquad
 \zeta\in\T\setminus U.
\]
Again putting $u=Dy$, we obtain
\[
 Q_X(\zeta,t)
 \geq
 \frac{m}{\norm{D}^2}I,
 \qquad
 \zeta\in\T\setminus U.
\]
Combining the two regions, there exists $c>0$ such that
\begin{equation}\label{eq:uniform-positive-pencil}
 Q_X(\zeta,t)\geq cI,
 \qquad
 \zeta\in\T.
\end{equation}

For $r<t$,
\[
 \begin{aligned}
 Q_X(\zeta,t)-Q_X(\zeta,r)
 &=
 \rho(t^2-r^2)I
 -2p(t-r)\Real(\overline\zeta X)\\
 &=
 (t-r)
 \Bigl(
 \rho(t+r)I-2p\Real(\overline\zeta X)
 \Bigr).
 \end{aligned}
\]
Hence
\[
 \sup_{\zeta\in\T}
 \norm{Q_X(\zeta,t)-Q_X(\zeta,r)}
 \leq
 (t-r)\bigl(\rho(t+r)+2p\norm{X}\bigr),
\]
which tends to zero as $r\uparrow t$. We may therefore choose
$r<t$ sufficiently close to $t$ so that
\[
 \sup_{\zeta\in\T}
 \norm{Q_X(\zeta,t)-Q_X(\zeta,r)}
 <c.
\]
Together with~\eqref{eq:uniform-positive-pencil}, this gives
\[
 Q_X(\zeta,r)\geq0,
 \qquad
 \zeta\in\T.
\]
Lemma~\ref{lem:positivity}{\rm(i)} then yields
\[
 w_\rho(X)\leq r<t,
\]
contradicting {\rm(i)}. Hence {\rm(ii)} holds.

\smallskip
\noindent\emph{{\rm(ii)}$\Rightarrow${\rm(i)}.}
Assume {\rm(ii)} and suppose, to the contrary, that
$
 w_\rho(X)<t.
$
Choose $r$ such that
\[
 w_\rho(X)<r<t.
\]
By Lemma~\ref{lem:positivity}{\rm(i)},
\[
 Q_X(\zeta,r)\geq0,
 \qquad
 \zeta\in\T.
\]
Equivalently,
\[
 Q_{X/r}(\zeta,1)\geq0,
 \qquad
 \zeta\in\T,
\]
so $X/r\in\mathcal C_\rho$.

Applying Theorem~\ref{thm:exact-DW} to $X/r$, we obtain
$
 \abs{\inner{Xx}{x}}\leq r
$
for every unit vector $x$. In particular,
\[
 \inner{\Real Xx}{x}
 =
 \Real\inner{Xx}{x}
 \leq
 \abs{\inner{Xx}{x}}
 \leq r,
\]
and hence
$
 \Real X\leq rI.
$

We now compare the boundary pencils at $r$ and $t$:
\[
 \begin{aligned}
 Q_X(1,t)-Q_X(1,r)
 &=
 (t-r)
 \bigl(\rho(t+r)I-2p\Real X\bigr)\\
 &\geq
 (t-r)
 \bigl(\rho(t+r)-2pr\bigr)I\\
 &=
 (t-r)(\rho t+qr)I.
 \end{aligned}
\]
Since $t>r>0$ and $q\geq0$, the last coefficient is strictly positive.
Thus, writing
\[
 \delta=(t-r)(\rho t+qr)>0,
\]
and using $Q_X(1,r)\geq0$, we obtain
\[
 Q_X(1,t)\geq\delta I.
\]

By~\eqref{eq:H-congruence}, for every $y\in\mathcal H$,
\[
 \begin{aligned}
 (1-s^2)\inner{H_{\rho,s}(C)y}{y}
 &=
 \inner{Q_X(1,t)Dy}{Dy}\\
 &\geq
 \delta\norm{Dy}^2\\
 &\geq
 \frac{\delta}{\norm{D^{-1}}^2}\norm{y}^2.
 \end{aligned}
\]
Therefore
\[
 H_{\rho,s}(C)
 \geq
 \frac{\delta}
 {(1-s^2)\norm{D^{-1}}^2}I.
\]
This contradicts
$0\in\sigma(H_{\rho,s}(C))$. Hence
\[
 w_\rho(X)=t=\ell_\rho(s),
\]
which proves {\rm(i)}.

\smallskip
\noindent\emph{Kernel characterization.}
Suppose first that
\[
 0\neq x\in\ker H_{\rho,s}(C),
\]
and normalize $x$. Put
\[
 z=\inner{Cx}{x},
 \qquad
 b=\norm{Cx}^2.
\]
Then~\eqref{eq:H-quadratic-form} gives
\[
 \varphi(z,b)=0.
\]
By the uniqueness of the zero of $\varphi$,
\[
 (z,b)=(z_{\rho,s},b_{\rho,s}).
\]
Lemma~\ref{lem:phase-support} shows that this point lies strictly on the
curved part of the upper boundary of $\mathfrak D_\rho$. Hence
Lemma~\ref{lem:curved-rigidity} provides a reducing two-dimensional
subspace $\mathcal M$ on which $C$ is unitarily equivalent to
$N_\rho$. Consequently,
\[
 C\simeq N_\rho\oplus C_0.
\]
Since $\mathcal M^\perp$ also reduces $C$, the compression lemma gives
$
 C_0\in\mathcal C_\rho.
$

Conversely, suppose that
\[
 C\simeq N_\rho\oplus C_0,
 \qquad
 C_0\in\mathcal C_\rho.
\]
The pair $(z_{\rho,s},b_{\rho,s})$ is realized by $N_\rho$, so there
exists a unit vector $x$ in the $N_\rho$-summand such that
\[
 \inner{Cx}{x}=z_{\rho,s},
 \qquad
 \norm{Cx}^2=b_{\rho,s}.
\]
By~\eqref{eq:H-quadratic-form},
\[
 \inner{H_{\rho,s}(C)x}{x}=0.
\]
Since $H_{\rho,s}(C)\geq0$,
\[
 0
 =
 \inner{H_{\rho,s}(C)x}{x}
 =
 \norm{H_{\rho,s}(C)^{1/2}x}^2,
\]
and therefore
\[
 H_{\rho,s}(C)x=0.
\]
Thus
\[
 \ker H_{\rho,s}(C)\neq\{0\}.
\]
This proves~\eqref{eq:kernel-reducing}.
\end{proof}

We now return directly to the original functional-calculus problem.

\begin{corollary}\label{cor:original-equality}
Let $1<\rho\leq2$, let $w_\rho(T)\leq1$, and let
$f\in A(\overline\D)$ satisfy
\[
 \norm{f}_{\infty,\D}\leq1.
\]
Assume $0<\abs{f(0)}<1$. Put $s=\abs{f(0)}$, choose $\omega\in\T$ such
that $\omega f(0)=s$, and define
\[
 h(z)=\frac{\omega f(z)-s}{1-s\omega f(z)}.
\]
Then
\[
 w_\rho(f(T))=\ell_\rho(s)
 \quad\Longleftrightarrow\quad
 0\in\sigma\bigl(H_{\rho,s}(h(T))\bigr).
\]
If $\mathcal H$ is finite-dimensional, these conditions are further
equivalent to
\[
 h(T)\simeq N_\rho\oplus C_0,
 \qquad
 C_0\in\mathcal C_\rho.
\]
\end{corollary}

\begin{proof}
For $z\in\overline\D$,
\[
 \abs{1-s\omega f(z)}
 \geq
 1-s\abs{f(z)}
 \geq
 1-s
 >
 0.
\]
Hence $h\in A(\overline\D)$. The defining disk automorphism gives
\[
 h(0)=0,
 \qquad
 \norm{h}_{\infty,\D}\leq1.
\]
Since $w_\rho(T)\leq1$, Lemma~\ref{lem:positivity}{\rm(i)} gives
$T\in\mathcal C_\rho$. Theorem~\ref{thm:zero-mapping} therefore yields
\[
 C:=h(T)\in\mathcal C_\rho.
\]

Moreover, from the definition of $h$,
\[
 \omega f=g_s\circ h,
\]
and hence, by the functional calculus,
\[
 \omega f(T)=g_s(C).
\]
Since $w_\rho$ is unimodularly invariant,
\[
 w_\rho(f(T))
 =
 w_\rho(\omega f(T))
 =
 w_\rho(g_s(C)).
\]
Theorem~\ref{thm:defect-equality} therefore gives
\[
 w_\rho(f(T))=\ell_\rho(s)
 \quad\Longleftrightarrow\quad
 0\in\sigma(H_{\rho,s}(C)),
\]
which is the first assertion.

Assume now that $\mathcal H$ is finite-dimensional. Then
\[
 0\in\sigma(H_{\rho,s}(C))
 \quad\Longleftrightarrow\quad
 \ker H_{\rho,s}(C)\neq\{0\}.
\]
By~\eqref{eq:kernel-reducing}, this is equivalent to
\[
 C\simeq N_\rho\oplus C_0,
 \qquad
 C_0\in\mathcal C_\rho.
\]
Since $C=h(T)$, the stated finite-dimensional characterization follows.
\end{proof}

The excluded endpoints have different equality structures. If
$|f(0)|=1$, then $f$ is constant. The case $f(0)=0$ reduces to the
equality problem for the normalized mapping theorem, while at
$\rho=1$ operator-norm equality has a different rigidity behavior.
\subsection*{Declarations}
\begin{itemize}
	\item \textbf{Availability of data and materials}: Data sharing is not applicable to this article.
	
	\item \textbf{Competing interests}: The authors declare that they have no competing interests.
	\item \textbf{Funding}: The authors received no funding for this work.
	\item \textbf{Authors' contributions}: All authors contributed equally and approved the final manuscript.
	\end{itemize}

\end{document}